\documentclass[3p]{elsarticle}

\usepackage{amsmath}
\usepackage{amssymb,amsbsy}
\usepackage{stackrel}

\newcommand{\R}{\mathbb R}
\newcommand{\g}[1]{\gamma_{#1}}
\newcommand{\blot}[1]{}

\renewcommand{\qed}{\hfill $\Box$}

\newtheorem{theo}{Theorem}
\newtheorem{lemma}[theo]{Lemma}

\begin{document}

\title{Point island dynamics under fixed rate deposition}
\tnotetext[1]{This paper is dedicated to the memory of Jack Carr}

\author[da]{D.~Allen}
\ead{damien.allen@strath.ac.uk}

\author[da]{M.~Grinfeld\corref{cor1}}
\ead{m.grinfeld@strath.ac.uk}

\author[rs1,rs2]{R.~Sasportes\fnref{fn1}}
\ead{rafael.sasportes@uab.pt}

\cortext[cor1]{Corresponding author}

\fntext[fn1]{Partially funded by FCT/Portugal through project
  RD0447/CAMGSD/2015.}

\address[da]{Department of Mathematics and Statistics, University of
  Strathclyde, Livingston Tower, 26 Richmond Street, Glasgow G1 1XH,
  UK}

\address[rs1]{Departamento de Ci\^encias e Tecnologia, Universidade
  Aberta, Lisboa, Portugal}
\address[rs2]{Centro de An\'alise Matem\'atica, Geometria e Sistemas
  Din\^amicas, Instituto Superior T\'ecnico, Universidade de Lisboa,
  Lisboa, Portugal}

\begin{abstract}
  We consider the dynamics of point islands during submonolayer
  deposition, in which the fragmentation of subcritical size islands
  is allowed. To understand asymptotics of solutions, we use methods
  of centre manifold theory, and for globalisation, we employ results
  from the theories of compartmental systems and of asymptotically
  autonomous dynamical systems. We also compare our results with those
  obtained by making the quasi-steady state assumption.
\end{abstract}  

\begin{keyword}
  submonolayer deposition \sep centre manifold theory \sep
  compartmental systems \sep asymptotically autonomous dynamical systems

\MSC[2010] 34A35 \sep 34E05 \sep 37C19 \sep 82D30
\end{keyword}
  
\maketitle

\section{Introduction}

Submonolayer deposition, a process in which atoms or molecules are
deposited onto a substrate, diffuse and form islands, is a
foundational technology in the creation of smart and nanomaterials
\cite{Mulheran2008}.  A mathematical theory of submonolayer deposition that
describes spatial distribution and the size statistics of the islands
is an important goal of research. At present there are many competing
models to describe the spatial distribution of islands; see, for
example \cite{Pimpinelli2007}, and the work that paper has led to.

Size distribution of islands is usually tackled by models that
disregard the spatial structure, and deal only with coagulation and
fragmentation of clusters composed of adatoms deposited onto a
surface. Such models lead to infinite systems of ordinary differential
equations (ODEs) for the various species, these are known as rate
equations; see, for example,  \cite{Einax2013, Einax2012}.

If furthermore one assumes that the structure of these clusters is
also disregarded, one deals with point islands, and then it makes
sense to assume that coagulation and fragmentation rates are not
size-dependent. Studies of this type of rate equations have been
initiated by da Costa {\em et al.}  \cite{daCosta2006}; see also
\cite{daCosta2016,Costin2013}, all of which are relevant to the
present work.
 
As in \cite{Costin2013} we further assume that there exists a {\em
  critical island size} $i$ such that islands (adatom clusters) of
size $j\geq n := i+1$ are immobile and can only grow by attachment of
single adatoms.

There is a number of possibilities how to model islands of size
$1<j \leq i$. The one considered in \cite{Costin2013} is that clusters
of size $1<j\leq i$ simply do not arise. There is one other physically
relevant possibility, i.e. that clusters of every size $1<j\leq i$ are
allowed to fragment (at some rate independent of the cluster size,
which is consistent with the point-island assumption). This
possibility has been considered formally in \cite{Blackman1991,
  MulheranBasham2008}. In this paper we consider this mechanism, using
centre manifold techniques \cite{Carr1981} and globalising the
results.

In \cite{daCosta2006} and in \cite{Costin2013} as well, it was
possible by a change of variables, to decouple the infinite system of
ODEs in a way that reduced its analysis to an analysis of a
two-dimensional system. In our case, the reduction is to $n=i+1$
equations, and the remarkable property of these equations is that the
complexity of the calculations is independent of $n$. Furthermore, it
appears that computations can be significantly simplified by making a
sweeping assumption that all the clusters of size $1<j \leq i$ are at
a quasi-steady state (the quasi-steady state assumption, QSSA). We
show that making this assumption results in the same leading term
behaviour as the centre manifold computation and emphasise the
differences between the two approaches.

The plan of the paper is as follows. In Section \ref{Eqns} we
introduce the governing equations, perform the finite-dimensional
reduction, and formulate an equivalence theorem between the reduced
$n$-dimensional system of equations and the original
infinite-dimensional one. In Section \ref{Glob} we discuss boundedness
and asymptotic behaviour of solutions to our equations. To obtain more
precise information about long-time asymptotics, in Section \ref{CM}
we perform a centre manifold analysis. The results in Section
\ref{Glob} imply that our asymptotics, derived by centre manifold
techniques, hold for any positive initial conditions. The
monomer asymptotics for large times are computed in Section
\ref{asymp} and are used there to discuss the consequences for the
asymptotic behaviour of islands of all sizes and to characterise the
similarity profile the solutions converge to.  In Section \ref{Equiv}
we compare our results to those obtained by making the QSSA, and
finally in Section \ref{Conc} we relate our results to those of
\cite{Blackman1991} and \cite{MulheranBasham2008} and draw
conclusions.

We also comment on the relation between the present paper and
\cite{daCosta2006}. The methods we use to obtain the asymptotics of
monomers and hence of the larger clusters (Lemmas \ref{c1a} and
\ref{cjli}) for all positive initial conditions, significantly extend the
methods used in \cite{daCosta2006}. It is in the way we use centre
manifold theory and globalise the results using the work of
\cite{Jacquez1993, Thieme1992} that the main novelty of the paper
lies. With the information contained in the above two lemmas, the
methods of \cite{daCosta2006} can be immediately used to discuss
similarity solutions (see Theorem \ref{simsol}); where proofs are
sufficiently similar to those in \cite{daCosta2006}, we either omit
them, or give only the gist as in the proof of Theorem \ref{exist}.

\section{Governing equations}\label{Eqns}

We consider a system containing clusters of any number $j \geq 1$ or
monomers. We assume that the following reactions occur:
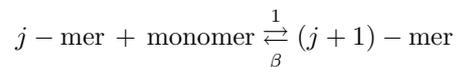
\[
  j-\hbox{mer}\,+\,\hbox{monomer} \stackrel[\beta]{1}{\rightleftarrows}
  (j+1)-\hbox{mer}
\]
for $1 \leq j < i$ and
\[
j-\hbox{mer}\,+\; \hbox{monomer}\stackrel[]{1}{\rightarrow} (j+1)-\hbox{mer}  
\]
if $j \geq i$. In other words, we allow clusters of size less than
$i+1$ to fragment at a rate $\beta>0$.  

If we set $\tilde{\alpha}$ to be the deposition rate, denote by $C_j(t)$ the
concentration of $j$-mers and use primes for differentiation with
respect to time $t$, the laws of mass kinetics give us the following
infinite system of ODEs:

\begin{equation}\label{odes1}
\begin{aligned}
C_1' & = \tilde{\alpha} -2 C_1^2 +2\beta C_2-C_1\sum_{k=2}^\infty C_k +\beta
\sum_{k=3}^i C_k,\\
C_j' & = C_1 C_{j-1}-C_1 C_{j} - \beta C_{j}+\beta
C_{j+1}, \quad 1<j<i,\\
C_j' &= C_1C_{j-1}-C_1C_j - \beta C_j, \quad j=i,\\
C_j' &= C_1C_{j-1}-C_1C_j, \quad j>i.
\end{aligned}
\end{equation}

It makes sense to scale the variables and the deposition rate to
remove $\beta$ from the equations. Thus scaling
$t\rightarrow T:= \beta t$, retaining primes for differentiation
with respect to the new time scale, setting $C_j(t)=\beta c_j(T)$ and
$\alpha = \tilde{\alpha}/\beta^2$, we obtain the system
\begin{equation}\label{odes2}
\begin{aligned}
c_1' & = \alpha -2 c_1^2 +2 c_2-c_1\sum_{k=2}^\infty c_k +\sum_{k=3}^i c_k,\\
c_j' & = c_1 c_{j-1}-c_1 c_{j} - c_{j}+ c_{j+1}, \quad 1<j<i,\\
c_j' &= c_1c_{j-1}-c_1c_j - c_j, \quad j=i,\\
c_j' &= c_1c_{j-1}-c_1c_j, \quad j>i.
\end{aligned}
\end{equation}

\section{Globalisation} \label{Glob}

In this section we consider the global dynamics of equations
(\ref{odes2}) satisfied by $c_j(T)$, $1 \leq j \leq i$, and
$v(T)=\alpha- c_1(T)\sum_{k=2}^\infty c_k(T)$, and establish that all
solutions of these equations with non-negative initial data approach
the origin. This will show that the flow on the centre manifold, as
given by Theorem \ref{asc1}, describes the asymptotics of every
non-negative solution of this system of equations.

For that purpose, it is more convenient to rewrite equations
(\ref{odes2}) formally as follows:
\begin{equation}\label{odes4}
\begin{aligned}
c_1' & = \alpha -2 c_1^2 +2 c_2-c_1\sum_{k=2}^i c_k +\sum_{k=3}^i
c_k -c_1y,\\
c_j' & = c_1 c_{j-1}-c_1 c_{j} - c_{j}+ c_{j+1}, \quad 1<j<i,\\
c_i' &= c_1c_{i-1}-c_1c_i - c_i, \\
y'   &= c_1 c_i,
\end{aligned}
\end{equation}
where we have put $y(T)= \sum_{k=i+1}^\infty c_k(T)$. 


First of all, we have 

\begin{theo}\label{exist}
If $\sum_{k=1}^\infty c_k(0) < \infty$, a solution of
\emph{(\ref{odes2})} for $j \geq 1$ is also a solution of
\emph{(\ref{odes4})}.
\end{theo}

{\em Proof:} The argument of the proof is similar to that of
\cite[Theorem 2.1]{daCosta2006}; we indicate the main steps. 

Let $(c_j)_{j=1}^{\infty}$ be a solution of (\ref{odes2}). To show
that this is also a solution of (\ref{odes4}) we must prove that
$\sum_{k=i+1}^{\infty}c_k$ converges to $y$ for all $T$. We change
time from $T$ to $\rho=\int_{0}^{T}c_1(s)\,{\rm d}s$. This change of
variable (also used in \cite[Theorem 2.1]{daCosta2006}) makes the
$c_j$ equations of (\ref{odes2}) linear in $c_j$ for $j>i$.  Keeping
primes for differentiation with respect to the new time variable
$\rho$ and letting $c_j(T):=\tilde{c}_j(\rho)$,
$y(T):= \tilde{y}(\rho)$, these equations become
\begin{equation}\label{eqrho}
  \tilde{c_j}'=\tilde{c}_{j-1}-\tilde{c}_j, \quad j> i,
  \hbox{ and } \tilde{y}'=\tilde{c}_i.
\end{equation}

This system of ODEs for $\tilde{c_j}$, $j> i$, can now be solved in
terms of $\tilde{c_i}$ recursively by variation of parameters starting at $j=i+1$, to give
\begin{equation}\label{eq:s3}
  \tilde{c}_j=e^{-\rho} \sum_{k=i+1}^{j}\frac{\rho^{j-k}}{(j-k)!}\tilde{c}_k(0) +\frac{1}{(j-(i+1))!} \int_{0}^{\rho} \tilde{c}_i (\rho-s) s^{j-(i+1)} e^{-s}\,{\rm d}s.
\end{equation}

Introducing  the generating function
\[
 F(\rho,z):= \sum_{n=i+1}^{\infty}\tilde{c}_nz^n,
\]
we can use (\ref{eq:s3}) to rewrite $F$ as
$F(\rho,z):=G(\rho,z)+H(\rho,z)$, where 
\[
G(\rho,z)=e^{-\rho} \sum_{n=i+1}^{\infty} \sum_{k=i+1}^{n}
\frac{\rho^{n-k} z^n}{(n-k)!}\tilde{c}_k(0),
\]
and
\[
  H(\rho,z)= \sum_{n=i+1}^{\infty}\frac{z^n}{(n-(i+1))!}
  \int_{0}^{\rho}\tilde{c}_i(\rho-s)s^{n-(i+1)}e^{-s}\,{\rm d}s.
\]

We now consider these two expressions separately.  For $G$ we obtain
\[
  G(\rho,z)=e^{-\rho(1-z)}\sum_{k=i+1}^{\infty}z^k \tilde{c}_k(0).
\]  
Since $ \sum_{k=1}^{\infty} c_k(0)<\infty$ by assumption, the above
series converges when $|z| \leq 1$, and we have
\[
G(\rho,z)=e^{-\rho(1-z)}(F(0,z)-\tilde{c}_i(0)z) \quad \hbox{ for }
\quad |z|\leq 1.
\]

For $H$, by interchanging the order of summation and integration, we
have that 
\[
H(\rho,z)=z^i \int_{0}^{\rho} \tilde{c}_1(s) e^{-(\rho-s)(1-z)}\,{\rm d}s.
\]
The expression for $F$ at $z=1$ now becomes
\begin{equation}\label{eqf}
F(\rho,1)= F(0,1)-\tilde{c}_i(0)+\int_{0}^{\rho}\tilde{c}_i(s)\,{\rm d}s.
\end{equation}
Hence, by differentiating with respect to
$\rho$, we see that
$F(\rho,1)$ given by (\ref{eqf}) satisfies the same differential
equation as $\tilde{y}$ in (\ref{eqrho}) which proves that $F(\rho, 1)
=\tilde{y}$. Thus in the $T$ variables
$\sum_{k=i+1}^{\infty}c_k$ converges to $y$.\qed

As a result of Theorem \ref{exist}, we can use finite-dimensional
techniques to discuss the dynamics of  $c_j(T)$, $1 \leq j \leq i$.  


We begin our analysis of long-time dynamics of (\ref{odes4}) by
considering the system without outflows through higher clusters, i.e.

\begin{equation}\label{eq:b1}
\begin{aligned}
c_1'&= \alpha -2c_1^2+2c_2-c_1\sum_{k=2}^{i} c_k +\sum_{k=3}^{i} c_k,\\
c_j'&=c_1c_{j-1}-c_1c_{j} -  c_{j}+ c_{j+1},   \; \;  1<j<i,\\
c_i'&=c_1c_{i-1}-  c_1 c_i- c_{i}.
\end{aligned}
\end{equation}

Let us show that the system (\ref{eq:b1}) is a compartmental system in
the sense of Jacquez and Simon \cite{Jacquez1993}. To that end, we introduce
some notation.

Let $I_1=\alpha$ and let $I_j=0$ for all $2 \leq j \leq i$.
Now put 
\begin{align*}
  F_{j1} & = c_1c_{j-1}, \; \; j=2, \ldots, i;\\
  F_{12} & = 2 c_2  \;\hbox{ and }\; F_{1j} = c_j ,\, j = 3, \ldots,
  i.
\end{align*}
For $k=j-1$, $2 \leq k \leq i-1$ put $F_{kj} = c_j$,
$F_{jk} = c_1 c_k$ and for $k=j+1$, $ 2 \leq j \leq i-1$, put
$F_{kj}=c_j$.  Finally, let $F_{0k}=0 $ if $k\neq 1,i$ and
$F_{0i}=F_{01}= c_1 c_i$, the only outflows from the system.

Then clearly for each $j = 1, \ldots, i$ we can write
\begin{equation}\label{comp}
  c_j'= \sum_{k \neq j}^i -F_{kj}+ F_{jk} +I_j-F_{0j},
\end{equation}  
where all the $F$s and $I$s are positive, which shows that
(\ref{eq:b1}) is a compartmental system in the sense of
\cite{Jacquez1993}.

Also note that
\begin{equation}\label{mc}
\frac{\partial F_{jk}}{\partial c_m}
\geq 0 \;\hbox{ for all }\; 1 \leq j,k,m \leq i,\; \; j \neq k.
\end{equation}

Hence we can use the theorem of Maeda, Kodama and Ohta \cite{Maeda1978}; see
also part (i) of Theorem 9 of \cite{Jacquez1993}:

\begin{theo}[\cite{Maeda1978}]\label{MKO}
  Given a compartmental system \emph{(\ref{comp})} with time-independent
  inputs $I_j$ that satisfies the monotonicity condition \emph{(\ref{mc})},
  every non-negative solution of the system is bounded iff the system has a
  positive rest point.
\end{theo}

Since it is not hard to compute  that the system (\ref{eq:b1}) admits the
unique positive equilibrium
\[
(c_1,c_2,\ldots,c_i)=(\alpha^{\frac{1}{i+1}},\alpha^{\frac{2}{i+1}}, \ldots, \alpha^{\frac{i}{i+1}}),
\]
we conclude using Theorem \ref{MKO} that all non-negative solutions of
(\ref{eq:b1}) are bounded.

Now we consider the first $i$ equations of the system
(\ref{odes4}). Since the system (\ref{odes4}) preserves
non-negativity, and $y(T)$ is a positive function, by comparison with
solutions of (\ref{eq:b1}) it follows that the $(c_1,\ldots, c_i)$
components of non-negative solutions of (\ref{odes4}) are bounded for
any positive initial condition.    

Now consider the dynamics of the last component of (\ref{odes4}),
$y(T)$. As it is monotone-increasing it can either converge to some
limit $l<\infty$ or it can go to infinity.   

Let us show that the first possibility cannot occur. If $y(T)$
converges to some limit $l<\infty$, we could use the theorem of Thieme
\cite[Theorem 4.2]{Thieme1992} on behaviour as $T \rightarrow \infty$
of asymptotically autonomous systems, combined with the fact that all
non-negative solutions of (\ref{odes4}) are bounded and the uniqueness
of the positive equilibrium, to conclude that the $\omega$-limit set
of every orbit of (\ref{odes4}) would be the same as that of the
system

\begin{equation}\label{odes5}
\begin{aligned}
c_1' & = \alpha -2 c_1^2 +2 c_2-c_1\sum_{k=2}^i c_k +\sum_{k=3}^i
c_k -c_1l,\\
c_j' & = c_1 c_{j-1}-c_1 c_{j} - c_{j}+ c_{j+1}, \; \; 1<j<i,\\
c_i' &= c_1c_{i-1}-c_1c_i - c_i.
\end{aligned}
\end{equation}

But if $y(T) \rightarrow l$ as $T \rightarrow \infty$, we must have
that either $c_1(T) \rightarrow 0$ or $c_i(T)\rightarrow 0$.
If we suppose, for example, that $c_i(T) \rightarrow 0$, we see from
the $c_i'$ equation of (\ref{odes5}) that either $c_1(T)$ or
$c_{i-1}(T)$ must go to zero. Continuing in this way, we see that all
$c_j(T)$ must go to zero as $T \rightarrow \infty$, but the origin is
not a rest point of the first $i$ equations of (\ref{odes5}). Hence we
conclude that $y(T) \rightarrow \infty$.  

Furthermore, since the positive orthant of $\R^{i+1}$ is invariant
under the flow of (\ref{odes4}), this  means that $c_1(T)
\rightarrow 0$ as $T \rightarrow \infty$.

Now, from the equations for $c_i\,,\: c_{i-1},\, \ldots,\, c_2$ it follows
consecutively that for all $2 \leq k \leq i$, $c_k(T) \rightarrow 0$
as $T \rightarrow \infty$, again using the same result of Thieme
\cite{Thieme1992} for asymptotically autonomous differential
equations. Applying these results to the equation for  $c_1$ in
(\ref{odes4}), we finally conclude that $c_1(T)y(T) \rightarrow \alpha$
as $T \rightarrow \infty$. If we now set 
\begin{equation}\label{defv}
v(T) := \alpha - c_1(T)\sum_{k=2}^\infty c_k(T), 
\end{equation} 
this is equivalent to saying that $v(T) \rightarrow 0$ as $T
\rightarrow \infty$.

We collect these results as a theorem:   

\begin{theo}\label{asall}
As $T  \rightarrow \infty$, for all non-negative initial data,
$c_k(T) \rightarrow 0$, $1\leq k \leq i$, and $v(T) \rightarrow 0$. 
\end{theo}  

To understand better the dynamics of $c_j(T)$ as $T \rightarrow
\infty$ for all $j \geq 1$, we first use centre manifold techniques to
understand the rate of approach of $c_j(T)$ to zero, $1 \leq i$ as $T
\rightarrow \infty$. 


\section{Centre manifold analysis} \label{CM}

The variable $v(T)$ defined by (\ref{defv}) satisfies the equation

\begin{equation}\label{odev}
\begin{aligned}  
  v' =& -\frac{1}{c_1} \Big[ c_1^4- c_1^2c_2+ \alpha v - v^2
- 2\alpha c_1^2 +2c_1^2v+2 \alpha c_2\\
&-2 c_2 v + \alpha \sum_{k=3}^i c_k -v \sum_{k=3}^i c_k \Big]. 
\end{aligned}
\end{equation}

In terms of $v$, the equation for $c_1$ becomes
\begin{equation}\label{c1}
c_1' = v - 2 c_1^2 +2 c_2 + \sum_{k=3}^i c_k.
\end{equation}

We now change time from $T$ to $\tau=\int_{0}^{T} \frac{1}{c_1(s)}\,
{\rm d}s$. This change of variable (also used in \cite[p. 377]{daCosta2006}
and \cite[system (3.3)]{daCosta2016}) is needed to desingularise the $v$
equation when $c_1=0$. Note that by the result of Theorem
\ref{asall}, $\tau \rightarrow \infty$ as $T \rightarrow \infty$.
Letting dots represent  derivatives with respect to $\tau$, the
$i+1$ ODEs for  $c_j$, $1 \leq j \leq i$, and for $v$  become
\begin{equation}\label{odes3}
\begin{aligned}
\dot c_1 & ~=~ c_1\left(v - 2 c_1^2 +2 c_2 + \sum_{k=3}^i c_k\right),\\
\dot c_j & ~=~ c_1(c_1 c_{j-1}-c_1 c_{j} - c_{j}+ c_{j+1}), \quad 1<j<i,\\
\dot c_i &~=~ c_1(c_1c_{i-1}-c_1c_i - c_i),
\end{aligned}
\end{equation}
and
\begin{equation}\label{odev1}
\begin{aligned} 
  \dot v =& -c_1^4- c_1^2c_2+ \alpha v - v^2
- 2\alpha c_1^2 +2c_1^2v\\
&~+2  \alpha c_2-2 c_2 v
+ \alpha \sum_{k=3}^i c_k -v \sum_{k=3}^i c_k. 
\end{aligned}
\end{equation}
Note that $0 \in \R^{i+1}$ is now a rest point of the system of equations
(\ref{odes3})--(\ref{odev1}). The object of interest is to establish
stability properties of this rest point and the way in which it is
approached.  


It is useful to make another change of variable. We set 
\[
w=v+ 2 c_2 +  \sum_{k=3}^{i} c_k.
\]

In the $(c_1,c_2,\dots, c_i,w)$ variables the equation (\ref{c1}) for
$c_1$ becomes conveniently
\begin{equation}\label{c1n}
\dot c_1= c_1 (w -2c_1^2),
\end{equation}
the equations for $c_2\,, \ldots,\, c_i$ remain as before in (\ref{odes3})
and obviously we have
\begin{equation}\label{odew}
\begin{aligned}
\dot w = &\dot{v} +2 \dot c_2 +\sum_{k=3}^{i} \dot c_k =-\alpha w +2
\alpha c_1^2 -2c_1 c_2+2 c_1^3 + c_1 c_3+4c_1^2c_2 - 2c_1^2 w\\ 
&~-2 c_2 w- c_1^2 c_i -w \sum_{k=3}^i c_k + 2c_1^2 \sum_{k=3}^i c_k-c_1^4 + w^2.
\end{aligned}
\end{equation}


Now we appeal to centre manifold theory \cite{Carr1981}. In the
language of that theory, for the equations (\ref{c1n}), (\ref{odes3}),
(\ref{odew}), the variables $c_j$, $ 1 \leq j \leq i$ are ``centre''
variables while $w$ is a ``stable'' variable, so that according to
centre manifold theory, in a neighbourhood of the origin in
$\R^{i+1}$, equations (\ref{c1n}), (\ref{odes3}) for $1<j \leq i$ and
(\ref{odew}) admit an $i$-dimensional centre manifold, $w=h(c_1,c_2,\,
\ldots, \,c_i)$. Furthermore, from Theorem \ref{asall} it follows that
the centre manifold attracts all solutions in a neighbourhood of the
origin in $\R^{i+1}$. 

On this centre manifold, the flow is given by 
\begin{equation}\label{cm1f}
\begin{aligned}
\dot c_1 &= c_1 (h(c_1,c_2,\ldots, c_i)-2c_1^2),\\
\dot c_j & = c_1(c_1 c_{j-1}-c_1 c_{j} - c_{j}+ c_{j+1}), \quad 1<j<i,\\
\dot c_i &= c_1(c_1c_{i-1}-c_1c_i - c_i).
\end{aligned}
\end{equation}

Remarkably, we can reparameterise time by going back to the 
variable $T$ to obtain on the centre manifold $w=h(c_1,c_2,\ldots, c_i)$
the equations
\begin{equation}\label{cm2f}
\begin{aligned}
c_1' &= h(c_1,c_2,\ldots, c_i)-2c_1^2,\\
c_j' &= c_1 c_{j-1}-c_1 c_{j} - c_{j}+ c_{j+1}, \quad 1<j<i,\\
c_i' &= c_1c_{i-1}-c_1c_i - c_i.
\end{aligned}
\end{equation}

Since by centre manifold theory the asymptotic expansion of
$h(c_1,c_2,\ldots, c_i)$ contains only quadratic terms and above in
$c_j$, $ j \geq 1$, the $i \times i$ Jacobian matrix $J(0)$ of
equations (\ref{cm2f}) around the origin in $\R^{i+1}$ has the
following structure:
\[
J(0)= \left[ \begin{array}{c|c}
              0 & 0 \\ \hline
              0 & A
        \end{array} \right],   
\]
with the first row being made of zeros, and the $(i-1) \times(i-1)$
bi-diagonal matrix $A$ having $-1$ on the main diagonal and $1$ in the
$(j,j+1)$ positions, $2 \leq j \leq i-1$.  It is easily seen that all
eigenvalues of $A$ are negative. Such structure of the Jacobian matrix
means that for the equations of the flow on the centre manifold $w=
h(c_1,c_2,\ldots, c_i)$, $c_j$, $2 \leq j \leq i\,$ are ``stable''
variables and $c_1$ is a ``centre'' variable, so that inside the
$i$-dimensional centre manifold there is another, one-dimensional,
centre manifold parameterised by $c_1$, i.e., a curve with components
$c_j= g_j(c_1)$, $1<j\leq i$. We will write
$g_w(c_1)=h(c_1,g_2(c_1),\ldots,g_i(c_1))$. Furthermore, we also know
by centre manifold theory that as $c_1 \rightarrow 0$,
\begin{equation}\label{cmcj}
g_j(c_1) \sim \sum_{k=2}^{\infty} \gamma_{j,k} c_1^k ,
\end{equation}
where we use $\sim$ to mean ``is asymptotic to as $c_1 \rightarrow
0$''. We also have 
\begin{equation}\label{cmw}
 g_w(c_1) \sim \sum_{k=2}^{\infty} \gamma_{w,k} c_1^k.
\end{equation} 

Hence (see \cite{Carr1981}) the flow on the one-dimensional centre
manifold is given by 
\begin{equation}\label{cm3f}
c_1' = g_w(c_1)-2c_1^2,
\end{equation}
and as the rest point at the origin of the one-dimensional ODE
(\ref{cm3f}) is asymptotically stable by Theorem \ref{asall}, the
one-dimensional centre manifold $(g_2(c_1), \ldots, g_i(c_1))$
attracts nearby solutions, so all (apart possibly from sets of zero
($i+1$)-dimensional Lebesgue measure) approach the origin along this
curve.

We have

\begin{theo}\label{asc1}
$c_1$ asymptotically satisfies the differential equation 
\[
c_1' \sim  \frac{1}{\alpha} \left( -c_1^{i+3}+
c_1^{i+4} - c_1^{2i+3}\right) + O(c_1^{2i+4}).
\]
\end{theo}
{\em Proof:} Before we start the computation of the coefficients
$\g{j,k}$ and $\g{w,k}$, let us indicate the flow of logic. The
equations we are dealing with, (\ref{cm2f}) and (\ref{odew}), have a
very special structure that we are going to exploit.

On the centre manifold, the equations determining $g_j(c_1)$, ($2 \leq j \leq i$), have the form

\begin{align*}
\frac{{\rm d}g_2}{{\rm d}c_1}(c_1) & (g_w(c_1)-2c_1^2)= c_1^2 - c_1 g_2(c_1)+c_1 g_3(c_1) -
g_2(c_1)+g_3(c_1),\\
\frac{{\rm d}g_j}{{\rm d}c_1}(c_1) & (g_w(c_1)-2c_1^2)= c_1 g_{j-1}(c_1)-c_1 g_j(c_1) -
g_j(c_1)+ g_{j+1}(c_1), \; 2<j <i,\\
\frac{{\rm d}g_i}{{\rm d}c_1}(c_1)& (g_w(c_1)-2c_1^2)= c_1 g_{i-1}(c_1)-c_1 g_i(c_1) -
g_i(c_1),
\end{align*}

to which we, denoting the right-hand side of  (\ref{odew}) by
$F(c_1,c_2, \ldots, c_i,w)$, add the equation
\[
\frac{{\rm d}g_w}{{\rm d}c_1}(c_1) (g_w(c_1)-2c_1^2)= F(c_1, g_2(c_1), \ldots,
g_i(c_1),g_w(c_1)).
\]
Now, using the expansions (\ref{cmcj}) and (\ref{cmw}), we have a
system of equations from which we can, in theory, find as many of the
coefficients $\gamma_{j,k}$ and $\gamma_{w,k}$ as we wish. The order
of the computation is as follows:

By inspection, one can immediately determine $\g{w,2}$, then
consecutively $\g{i,2}$, $\g{i-1,2}$ and all the way to
$\g{2,2}$. Once this is done, we can find $\g{w,3}$ and proceed in
this way to find as many terms of the expansion as required. See 
Appendix A for the MAPLE code to do the computation for $i=5$.

Following the algorithm, we find that for all $j$, $2\leq j\leq i$,
$g_j(c_1) = O(c_1^j)$ and $\g{j,j}=1$. The final result of this
computation is that the functions $g_j(c_1)$, $2\leq j< i$,
$g_i(c_1)$, and $g_w(c_1)$ have the following asymptotic expansions as
$c_1 \rightarrow 0$:
\begin{equation}\label{gj}
  g_j(c_1) \sim c_1^j -c_1^{i+1}+c_1^{i+j} + O(c_1^{i+j+2}), \; \;
  g_i(c_1) \sim c_1^i-c_1^{i+1} + c_1^{2i} + O(c_1^{2i+1}),
\end{equation}
and 
\[
g_w(c_1) \sim 2c_1^2 + \frac{1}{\alpha} \left(-c_1^{i+3}+
c_1^{i+4} - c_1^{2i+3} \right)+ O(c_1^{2i+4}). 
\]
From these representations, Theorem \ref{asc1} follows 
immediately. \qed

Note that beyond terms of $O(c_1^{i+j})$ the interplay among
$g_j(c_1)$, $1<j\leq i$, and $g_w(c_1)$ becomes complex, and that the
later coefficients of these functions depend on $\alpha$. Computations
using the MAPLE code in  Appendix A indicate that the radius of
convergence of the expansions is $0$ for all $\alpha>0$.


\section{Asymptotics of solutions}\label{asymp}

Armed with Theorem \ref{asc1}, which holds for any non-negative
solution of (\ref{odes2}) by the globalisation results of Section
\ref{Glob}, we can discuss asymptotics of solutions of (\ref{odes1}),
using the methods of \cite{daCosta2016,daCosta2006}, which were also
used in \cite{Costin2013}. As proofs are similar to those used in the
above papers, we indicate only the main ideas.  Further terms in the
expansions in this section can be computed using the machinery of
\cite{daCosta2016}; here we only determine the leading terms, denoting
higher order terms by ``h.o.t.''. Going back to our original variables
$C_j(t)$ to exhibit the complicated dependence of the results on
$\beta$, from Theorem \ref{asc1} we have the following statement:

\begin{lemma}\label{c1a}
As $t \rightarrow \infty$, the asymptotics of $C_1(t)$ are given by 
\[
C_1(t) \sim  \bigg(\frac{\tilde{\alpha} \beta^{i-1}
}{(i+2)t}\bigg)^\frac{1}{i+2} + \emph{\hbox{h.o.t.}}
\]
\end{lemma}

Note that if we set $\beta=1$ in the equation above we obtain the same
result as in \cite{Costin2013}. Already at the level of $C_1(t)$ one sees
that the influence  of the fragmentation rate $\beta$ is not
intuitive.

Once we know the asymptotics of $C_1(t)$ from Lemma \ref{c1a}, the
asymptotics of $C_j(t)$ when $ 1 \leq j \leq i$ follow from (\ref{gj}).

\begin{lemma}\label{cjli}
 For $ 1 < j \leq i$, the asymptotics of $C_j(t)$ as $t \rightarrow
 \infty$ are given by
\[
  C_j (t)\sim \bigg(\frac{\tilde{\alpha} \beta^{\frac{i-3j+2}{j}}
  }{(i+2)t}\bigg)^{\frac{j}{i+2}}+\emph{\hbox{h.o.t.}}
\]
\end{lemma} 

Hence we are now in a position to express the asymptotics of $C_j(t)$
when $j>i$ by solving linear non-homogeneous ODEs using the same change of
variable as in the proof of Theorem \ref{exist}. 

\begin{lemma}\label{cjgi}
For $ j > i$, the asymptotics of $C_j(t)$  as $t \rightarrow
 \infty$ are given by 
\[
C_j (t)\sim \bigg(\frac{\tilde{\alpha} \beta^{\frac{-2i+2}{i}}
}{(i+2)t}\bigg)^{\frac{i}{i+2}}+\emph{\hbox{h.o.t.}}
\]
\end{lemma}

From this information we have the equivalent of \cite[Theorem
5.1]{daCosta2006} and, to which it is more directly comparable,
\cite[Theorem 6]{Costin2013} concerning similarity solutions of
(\ref{odes1}). These references should be consulted for the required
computations. To formulate the theorem, we first compute the
asymptotics of the average cluster size $\langle j \rangle$ using the
information in Lemmas \ref{c1a}--\ref{cjgi}:
\[
\langle j \rangle   = \frac{\sum_{j=1}^\infty j
  C_j(t)}{\sum_{j=1}^\infty C_j(t)} \sim
\left( \frac{\tilde{\alpha} \beta ^{i-1}}{i+2} \right)^{\frac{1}{i+2}}
t^{\frac{i+1}{i+2}} + \hbox{h.o.t}.
\]
Next, we define the function $\Psi$ by
\[
  \Psi (r) = \begin{cases}
    (1-r)^{-\frac{i}{i+1}} & if $r< 1$,\\
    0 & \hbox{otherwise}.
  \end{cases}
\]  

Finally, define the similarity variable $\eta$ by
\[
  \eta = \frac{(i+1)\beta^{-\frac{i+1}{i+2}}}{i+2}\frac{j}{\langle
          j \rangle}.
\]
Then we have that the solutions of (\ref{odes1}) converge to a
(discontinuous) similarity profile:

\begin{theo}\label{simsol}
  \(
    C_j(t) =\langle j \rangle^{-\frac{i}{i+1}} \Psi
      \left( \eta \right)\)  as $t \rightarrow \infty$.
\end{theo}  

The profile obtained in this theorem can be further analysed by the
methods of \cite[Section 6]{daCosta2006}.

\section{Quasi-steady state assumption} \label{Equiv}

In this section we would like to investigate whether the asymptotics
of solutions obtained in Section \ref{asymp} based on the centre
manifold analysis of Section \ref{CM} can be recovered more easily by
combining centre manifold reasoning with a technique that is often
used in the engineering community, the quasi-steady state
approximation (QSSA; see \cite{Goussis2012, Pantea2014,Segel1989}). As
in the famous example from enzyme kinetics due to Segel and Slemrod
\cite{Segel1989}, we show that QSSA correctly captures the leading
term asymptotics, though of course there will be differences in higher
order terms.

We restart with equations (\ref{odes2}), but now we immediately make
the QSSA assumption that $c_j' =0$ for $1< j \leq i$. We solve
the $i$ algebraic equations for $c_j$, $1< j \leq i$,
in terms of $c_1$. This clearly can be done consecutively, by starting
with the $c_{i}$ equation and solving it in terms of $c_1$ and 
$c_{i-1}$, substituting the expression we get for $c_{i}$
into the $c_{i-1}$ equation and continuing in this way, till $c_2$ has
been solved in terms of $c_1$, after which we back-substitute.

This procedure gives us that under the QSSA assumption
\[
c_j = \frac{\sum_{k=1}^{i-j+1} c_1^{k+j-1}}{\sum_{k=1}^{i}
  c_1^{k-1}},\: j = 2, \ldots, i. 
\]
Note that these are global objects, defined for all values of $c_1
>0$ unlike the centre manifold expansions (\ref{gj}) which have zero
radius of convergence.   We will need the MacLaurin series expansions
of these objects, 
\begin{equation}\label{qssa}
  c_j = c_1^j + \sum_{k=1}^n \left(-c_1^{ki+1} + c_1^{ki+j} \right) +
  O\left(c_1^{(n+1)i+1}\right).
\end{equation}

Now we can go back to the equation for $\dot w$ (\ref{odew}), write
$w=g_w(c_1)$, remember that by the centre manifold theorem $g_w$
contains terms that are at least quadratic in $c_1$, and substituting
instead of $c_j$, $1< j \leq i$, the expressions from (\ref{qssa}),
obtain that
\[
  g_w(c_1) \sim 2 c_1^2 + \frac{1}{\alpha} \sum_{k=1}^n
  \left(-c_1^{ki+3} + c_1^{ki+4} \right) + O \left(c_1^{(n+1)i+3}\right).
\]
In the series above we can take $n$ as large as we wish. Hence under
the QSSA assumption, setting $h(c_1,\ldots,\,c_i)=g_w(c_1)$ in the
first equation of (\ref{cm2f}), the dynamics of $c_1$ is governed by
the equation
\[
c_1' \sim \frac{1}{\alpha} \left( - c_1^{i+3} + c_1^{i+4} - c_1^{2i+3}
  + c_1^{2i+4} \right) + O \left(c_1^{3i+3}\right),
\]
which by inspection yields the same first three terms as the centre
manifold computation of Theorem \ref{asc1} for a fraction of the
effort.

\section{Conclusions} \label{Conc}  

In this paper we complemented the analysis of \cite{Costin2013} by
considering a more realistic dynamics of nucleating point islands with
critical island size $i$ by allowing subcritical islands of size
$2 \leq j \leq i-1$ to form and fragment. The mathematics of this new
system of equations is more challenging than the fundamentally
2-dimensional system considered in \cite{Costin2013} and we had to use both
centre manifold techniques and a sophisticated globalisation argument
using ideas from theories of compartmental systems and of asymptotically
autonomous differential equations; the globalisation methods used in
this paper are in our opinion more elegant than the ``brute-force''
asymptotics in \cite{Costin2013}.  

Our asymptotic results in Section \ref{asymp} are consistent with the
leading term asymptotics for $c_1(t)$ of \cite{MulheranBasham2008}
(see our Lemma \ref{c1a}) and for $c_j(t)$, ($1 \leq j \leq i$), of
\cite{Blackman1991} (see Lemma \ref{cjli}), as well as with the
conjecture in \cite{MulheranBasham2008} about the behaviour of
$c_j(t)$, $j > i$ (see Lemma \ref{cjgi}). Of course our methods are
not restricted to the computation of leading terms of the asymptotics.

\appendix 

\section{Computations of Theorem \ref{asc1}}

In this Appendix we supply the code implementing the computations
described in the proof of Theorem \ref{asc1}. We compute $15$ terms in
the expansion of the one-dimensional centre manifold with components
given by (\ref{cmcj}) of a system with $i=5$.  

\begin{verbatim} 
n:=15:
\end{verbatim}

First of all we set up the equations:

\begin{verbatim}

eqc1 := a-2*c1^2-c1*z+2*b*c2+b*c3+b*c4+b*c5:
eqc2 := c1^2-c1*c2-b*c2+b*c3:
eqc3 := c1*c2-c1*c3-b*c3+b*c4:
eqc4 := c1*c3-c1*c4-b*c4+b*c5:
eqc5 := c1*c4-c1*c5-b*c5:

eqz := c1^2-b*c2:

z := (a-v)/c1:

eqv := -eqc1*z-c1*eqz:

eqc1s := eqc1*c1:
eqc2s := eqc2*c1:
eqc3s := eqc3*c1:
eqc4s := eqc4*c1:
eqc5s := eqc5*c1:

eqvs := simplify(eqv*c1):

v := w-2*b*c2-b*c3-b*c4-b*c5:

eqws := simplify(eqvs+2*b*eqc2s+b*eqc3s+b*eqc4s+b*eqc5s):

eqws := expand(eqws):

\end{verbatim}

Now we use the expansions (\ref{cmcj}) and (\ref{cmw}) and remove all
the higher order terms that are not needed in the computation to save
time:

\begin{verbatim} 

c2 := sum('g2||j*c1^j','j'=2..n): 
c3 := sum('g3||j*c1^j','j'=2..n): 
c4:= sum('g4||j*c1^j','j'=2..n): 
c5 := sum('g5||j*c1^j','j'=2..n):
w := sum('gw||j*c1^j','j'=2..n):

aw := collect(simplify(eqc1s*diff(w,c1)-eqws),c1):
ac2 := collect(simplify((w-2*c1^2)*diff(c2,c1)-eqc2),c1):
ac3 := collect(simplify((w-2*c1^2)*diff(c3,c1)-eqc3),c1):
ac4 := collect(simplify((w-2*c1^2)*diff(c4,c1)-eqc4),c1):
ac5 := collect(simplify((w-2*c1^2)*diff(c5,c1)-eqc5),c1):
aw:= convert(taylor(aw,c1=0,n+1),polynom):

ac2:= convert(taylor(ac2,c1=0,n+1),polynom):
ac3:= convert(taylor(ac3,c1=0,n+1),polynom):
ac4:= convert(taylor(ac4,c1=0,n+1),polynom):
ac5:= convert(taylor(ac5,c1=0,n+1),polynom):

\end{verbatim}

Finally, we compute the coefficients of the expansion in the order
indicated in the proof of Theorem \ref{asc1}. 

\begin{verbatim} 

for k from 2 to n do 
     gw||k:= solve(coeff(aw,c1,k),gw||k):
     g5||k:= solve(coeff(ac5,c1,k),g5||k):
     g4||k:= solve(coeff(ac4,c1,k),g4||k):
     g3||k:= solve(coeff(ac3,c1,k),g3||k):
     g2||k:= solve(coeff(ac2,c1,k),g2||k):
od:

\end{verbatim}

Now we print out the asymptotic ODE equation for $c_1$:

\begin{verbatim}
odec1 := w-2*c1^2;
\end{verbatim} 

The result is 
\[
c_1' \sim 
- \frac{c_1^8}{\alpha \beta^4} + \frac{c_1^9}{\alpha\beta^5} 
- \frac{c_1^{13}}{\alpha\beta^9} + 
\frac{(30\beta^2+\alpha)c_1^{14}}{\alpha^2\beta^{10}}  - 
\frac{80c_1^{15}}{\alpha\beta^9}+ O\left(c_1^{16}\right). 
\] 

\medskip

{\bf Acknowledgment:} Part of this work was done during MG's visit to
the Department of Mathematics of the University of Aveiro, whose
hospitality is gratefully acknowledged. RS acknowledges the support of
the Edinburgh Mathematical Society to Glasgow where the work was
initiated. The comments of an anonymous referee are gratefully
acknowledged; they have led to a substantial improvement in the
readability of the paper. 

\medskip 
\section*{References}


\end{document}